\documentclass{amsart}

\usepackage[dvips]{color}

\usepackage{amsmath}
\usepackage{amssymb}

\theoremstyle{plain}   %% This is the default, anyway
\begingroup % Confine the \theorembodyfont command
\newtheorem{theorem}{Theorem}   % Numbered within each section
\newtheorem{corollary}[theorem]{Corollary}     % Numbered along with thm
\newtheorem{lemma}[theorem]{Lemma}         % Numbered along with thm
\endgroup

\theoremstyle{definition}
\newtheorem{definition}[theorem]{Definition}   % Numbered along with thm

\theoremstyle{remark}
\newtheorem{remark}[theorem]{Remark}        % Numbered along with thm
\newcommand{\Rn}{{\mathbb R}^{n}}

\newcommand{\R}{{\mathbb R}}
\newcommand{\N}{{\mathbb N}}

\newcommand{\Q}{\mathbb Q}
\newcommand{\e}{\varepsilon}

\newcommand{\lin}{\operatorname{span}}

\newcommand{\vp}{\varphi}

\newcommand{\Lip}{\operatorname{Lip}}

\newcommand{\interior}{\operatorname{int}}

\newcommand{\mcA}{{\mathcal{A}}}
\newcommand{\mcC}{{\mathcal{C}}}

\begin{document}
\title{On G\^ateaux differentiability of
pointwise Lipschitz mappings}



\author{Jakub Duda}
\thanks{The author was supported in part by ISF}

\email{duda@karlin.mff.cuni.cz}

\address{
Department of Mathematics\\
Weizmann Institute of Science\\
Rehovot 76100\\
Israel}

\curraddr{
Charles University,
Department of Mathematical Analysis,
Sokolovsk\'a 83,
186 75 Praha 8-Karl\'\i n,
Czech Republic}

\date{July 26, 2006}
\begin{abstract} 
We prove that for every function $f:X\to Y$,
where $X$ is a separable Banach space and $Y$ is a Banach space
with RNP, there exists a set $A\in\tilde\mcA$ such that $f$ is
G\^ateaux differentiable at all $x\in S(f)\setminus A$, where
$S(f)$ is the set of points where $f$ is pointwise-Lipschitz.
This improves a result of Bongiorno. As a corollary,
we obtain that every $K$-monotone function on a separable Banach space
is Hadamard differentiable outside of a set belonging to $\tilde\mcC$;
this improves a result due to Borwein and Wang. 
Another corollary is that if $X$ is Asplund, $f:X\to\R$ cone monotone,
$g:X\to\R$ continuous convex, then there exists a point in $X$, where $f$ is Hadamard
differentiable and $g$ is Fr\'echet differentiable.
\end{abstract}

\maketitle

\section{Introduction}
The classical Rademacher theorem~\cite{R} concerning a.e.\ differentiability
of Lipschitz functions defined on $\Rn$ was
extended by Stepanoff to pointwise Lipschitz functions~\cite{S1,S2}.
D.~Bongiorno~\cite[Theorem~1]{B} proved a version for infinite-dimensional mappings;
namely, that for every $f:X\to Y$, where $X$ is a separable Banach space
and $Y$ is a Banach space with RNP, there exists an Aronszajn null
set $A\subset X$ (see e.g.\ \cite{BL} for the definition of Aronszajn null sets) 
such that $f$ is G\^ateaux
differentiable at all $x\in S(f)\setminus A$ (here, $S(f)$ is the
set of points where $f$ is pointwise-Lipschitz). 
This generalized results for Lipschitz functions obtained by
Aronszajn, Christensen, Mankiewicz, and Phelps; see e.g.\ \cite{BL}
for the definitions of various notions of null sets they used.
We prove a stronger version of infinite dimensional Stepanoff-like theorem,
which asserts that under the same assumptions as in~\cite[Theorem~1]{B}, the set $A$ can
be taken in the class $\tilde\mcA$ defined by Preiss and Zaj\'\i\v cek~\cite{PZ};
see Theorem~\ref{mainthm}.
By results of~\cite{PZ}, $\tilde\mcA$ is a strict subclass of Aronszajn null sets.
Recently, Zaj\'\i\v{c}ek~\cite{Znew} proved that the
sets in $\tilde\mcA$ (and even $\tilde\mcC$) are $\Gamma$-null, which is a notion of null sets
due to Lindenstrauss and Preiss~\cite{LP} (here, a definition and basic
properties of this notion can be found).
Thus, Theorem~\ref{mainthm} has the following corollary:
if $X$ is a Banach space with separable dual (i.e.\ an Asplund space),
and $Y$ is a Banach space with RNP, $f:X\to Y$ is pointwise-Lipschitz at all $x\in X\setminus A$
where $A\in\tilde C$,
$g:X\to\R$ is continuous convex, then there exists $x\in X$ such that $f$ is G\^ateaux differentiable
at~$x$ and $g$ is Fr\'echet differentiable at~$x$.
In some sense, our proof of Theorem~\ref{mainthm} is simpler
than the proof of~\cite[Theorem~1]{B}; some of the (rather cumbersome) measurability
considerations from~\cite{B} are replaced by Lemma~\ref{Oborel}
and the construction of a total set from~\cite{B} is
replaced by the Lipschitz property of certain restrictions of the given mapping.
In the proof, we use several ideas from~\cite{PZ}.	
\par
Let $X$ be a Banach space and $\emptyset\neq K\subset X$ be a cone.
Following~\cite{BBL}, we say that $f:X\to\R$ is {\em $K$-monotone} 
provided $f$ or $-f$ is $K$-increasing (we say that $f:X\to\R$ is
{\em $K$-increasing} provided $x\leq_K y$ implies $f(x)\leq f(y)$ whenever
$x,y\in X$; here, $x\leq_K y$ means $y-x\in K$).
Borwein, Burke and Lewis~\cite{BBL} proved that every $K$-monotone
$f:X\to\R$ is G\^ateaux differentiable outside of a Haar null set
(see~\cite{BL} for definition) provided $X$ is separable and $K$ is
closed convex with $\interior(K)\neq\emptyset$.
This was strengthened by Borwein and Wang~\cite{BW} who showed that
``Haar null'' can be replaced by ``Aronszajn null''.
In section~\ref{monotonesection}, as a corollary to Theorem~\ref{mainthm},
we obtain that an analogous result holds if we replace ``Haar null''
by the class $\tilde\mcC$ defined by Preiss and Zaj\'\i\v{c}ek~\cite{PZ};
see Theorem~\ref{conethm} for details. 
The class $\tilde\mcC$ is a strict subclass of Aronszajn null sets
(see~\cite[p.\ 19]{PZ}) and thus our result improves the result due
to Borwein and Wang.
\cite[Proposition~16(iv)]{BW} shows that instead of ``G\^ateaux
differentiable'' we can write ``Hadamard differentiable'' (see Corollary~\ref{Hadcor}).
Our result has another interesting corollary; namely, if $X$ has a separable
dual (i.e.\ $X$ is an Asplund space), 
$f:X\to\R$ is $K$-monotone, $g:X\to\R$ is continuous convex, then
there exists $x\in X$ such that $f$ is Hadamard differentiable at~$x$,
and $g$ is Fr\'echet differentiable at $x$ (see Corollary~\ref{corgamma}).
This does not follow from
the results of Borwein and Wang since Aronszajn null sets and $\Gamma$-null
sets are incomparable.
It seems to be a difficult
open problem whether $\tilde\mcC=\tilde\mcA$ (see~\cite{PZ}). If this were true,
then our theorem would also hold with $\tilde\mcA$ in place of $\tilde\mcC$.
Thus, it remains open, whether we can replace $\tilde\mcC$ by $\tilde\mcA$
in Theorem~\ref{conethm} and Corollary~\ref{Hadcor}.
Going in another direction, the author~\cite{Dcone} proved some
results about a.e.\ differentiability of vector-valued cone monotone mappings.
\par
The current paper is organized as follows. Section~\ref{prelims}
contains basic definitions and facts. Section~\ref{auxiliary}
contains auxiliary results. Section~\ref{mainsec} contains
the proofs of the main Theorem~\ref{mainthm}, and Corollary~\ref{maincor}.
Section~\ref{monotonesection} contains the proofs of Theorem~\ref{conethm},
and Corollaries~\ref{Hadcor} and~\ref{corgamma}.

\section{Preliminaries}\label{prelims} 

All Banach spaces are assumed to be real.
By $\lambda$ we will denote the Lebesgue measure on $\R$.
Let $X$ be a Banach space. By $B(x,r)$ we will denote
the open ball with center $x\in X$ and radius $r>0$,
and by $S_X$ we denote $\{x\in X:\|x\|=1\}$.
If $M\subset X$, then by $d_M(x):=\inf\{\|y-x\|:y\in M\}$
we denote the distance from $x\in X$ to~$M$.
\par
Let $X,Y$ be Banach spaces.
We say that $f:X\to Y$ is {\em pointwise Lipschitz at $x\in X$}, provided
$\limsup_{y\to x}\frac{\|f(x)-f(y)\|}{\|x-y\|}<\infty$.
By $S(f)$, we will denote the set of points of~$X$ where~$f$ is pointwise Lipschitz.
By $\Lip(f)$ we will denote the usual Lipschitz constant of~$f$.
\par
In the following, let $X$ be a Banach space. 
If $f$ is a mapping from $X$ to a Banach space $Y$ and $x,v\in X$,
then we consider the directional derivative $f'(x,v)$ defined by
\begin{equation}\label{gateq} 
f'(x,v)=\lim_{t\to0}\frac{f(x+tv)-f(x)}{t}.
\end{equation}
If $x\in X$, $f'(x,v)$ exists for all $v\in X$, and
$T(v):=f'(x,v)$ is a bounded linear operator from~$X$ to~$Y$, then
we say that {\em $f$ is G\^ateaux differentiable} at~$x$.
If $f$ is G\^ateaux differentiable at $x$ and the limit in~\eqref{gateq} is
uniform in $\|v\|=1$, then we say that $f$ is {\em Fr\'echet differentiable} at~$x$.
If $f$ is G\^ateaux differentiable at~$x$, and the limit in~\eqref{gateq}
is uniform with respect to norm-compact sets, then we say that
$f$ is {\em Hadamard differentiable} at $x$.
\par
We will need the following notion of ``smallness'' of sets
in Banach spaces from~\cite{PZ}.
\begin{definition}
Let $X$ be a Banach space, $M\subset X$, $a\in X$. Then
we say that
\begin{enumerate}
\item $M$ is {\em porous} at $a$ if there exists $c>0$ such that 
for each $\e>0$ there exist $b\in X$ and $r>0$ such that
$\|a-b\|<\e$, $M\cap B(b,r)=\emptyset$, and $r>c\|a-b\|$.
\item $M$ is {\em porous at $a$ in direction $v$} if the
$b\in X$ from~(i) verifying the porosity of $M$ at $a$ can
be always found in the form $b=a+tv$, where $t\geq0$. We say
that $M$ is {\em directionally porous at~$a$} if there
exists $v\in X$ such that $M$ is porous at~$a$ in direction~$v$.
\item $M$ is {\em directionally porous} if $M$ is directionally
porous at each of its points.
\item $M$ is {\em $\sigma$-directionally porous} if it is
a countable union of directionally porous sets.
\end{enumerate}
\end{definition}
For a recent survey of properties of negligible sets, see~\cite{Zajsur}.
We will also need the following notion of ``null'' sets in a Banach space.
It was defined in~\cite{PZ}.

\begin{definition} 
Let $X$ be a separable Banach space and $0\neq v\in X$.
Then $\tilde\mcA(v,\e)$ is the system of all Borel sets $B\subset X$
such that $\{t:\vp(t)\in B\}$ is Lebesgue null whenever $\vp:\R\to X$ is such
that the function $t\to\vp(t)-tv$ has Lipschitz constant at most~$\e$, 
and $\tilde\mcA(v)$ is the system of all sets $B$ such that
$B=\bigcup_{k=1}^\infty B_k$, where $B_k\in\tilde\mcA(v,\e_k)$ for some $\e_k>0$.
\par
We define $\tilde\mcA$ (resp.\ $\tilde\mcC$) as the system of those $B\subset X$ that can be,
for every given complete\footnote{We say that $(v_n)_n\subset X\setminus\{0\}$ is
a {\em complete} sequence provided $\overline{\lin}(v_n)=X$.} sequence $(v_n)_n$ in $X$
(resp.\ for some sequence $(v_n)_n$ in $X$), written
as $B=\bigcup_{n=1}^{\infty} B_n$, where each $B_n$ belongs
to $\tilde\mcA(v_n)$.
\end{definition}

The following simple lemma shows that every directionally
porous set is contained in a set from $\tilde\mcA$. As a corollary, we have
the same result for $\sigma$-directionally porous sets.

\begin{lemma}\label{borellem}
Let $X$ be a separable Banach space, and $A\subset X$ be directionally porous.
Then there exists a set $\hat A\in\tilde\mcA$ such that $A\subset\hat A$.
\end{lemma}

\begin{proof} 
This follows from the proof of~\cite[Theorem~10]{PZ};
see also~\cite[Remark~6]{PZ}.
\end{proof}
\par
The following simple lemma is proved in~\cite{B}:
\begin{lemma}[\cite{B}, Lemma~1]\label{bong1lem}
Given $f:X\to Y$ and $L,\delta>0$, let $S$ be the set of all points $x\in X$
such that 
$\|f(x+h)-f(x)\|\leq L\|h\|$ whenever $\|h\|<\delta$.
Then $S$ is a closed set.
\end{lemma}

\section{Auxiliary results}\label{auxiliary}

The following is an extension of~\cite[Lemma~3]{BW} to vector-valued
setting.

\begin{lemma}\label{likebwlem}
Let $X,Y$ be Banach spaces, $f:X\to Y$.
Fix $v_1,v_2\in X$, $k,l,m\in\N$, and $y,z\in Y$. Then the set
$A(k,l,m,y,z)$ of all $x\in X$ verifying
\begin{enumerate}
\item 
$\big\|\frac{f(x+tu)-f(x)}{t}-y\big\|<\frac{1}{l}$
for $\|u-v_1\|<1/m$
and $0<t<1/k$;
\item 
$\big\|\frac{f(x+tu)-f(x)}{t}-z\big\|<\frac{1}{l}$
for $\|u-v_2\|<1/m$ and $0<t<1/k$; and
\item 
$\big\|\frac{f(x+s(v_1+v_2))-f(x)}{s}-(y+z)\big\|>\frac{3}{l}$
occurs for arbitrarily small $s>0$,
\end{enumerate}
is directionally porous in $X$.
\end{lemma}

\begin{proof} Let $x\in A(k,l,m,y,z)$. Choose $0<s<1/k$ such that the inequality in~(iii) holds.
We claim that 
$B\big(x+sv_1,\frac{s}{m}\big)\cap A(k,l,m,y,z)=\emptyset$.
\par
Indeed, for~$\|h\|<\frac{1}{m}$, if $x+s(v_1+h)$ satisfies~(ii), we have
\begin{equation}\label{bw1eq}
\bigg\|\frac{f(x+s(v_1+h)+su)-f(x+s(v_1+h))}{s}-z\bigg\|<\frac{1}{l},
\end{equation}
for $\|u-v_2\|<\frac{1}{m}$.
By~(i) we get
\begin{equation}\label{bw2eq}
\bigg\|\frac{f(x+s(v_1+h))-f(x)}{s}-y\bigg\|<\frac{1}{l}.
\end{equation}
By the triangle inequality,~\eqref{bw1eq}, and~\eqref{bw2eq} we get
\[
\bigg\|\frac{f(x+s(v_1+h)+su)-f(x)}{s}-(y+z)\bigg\|<\frac{2}{l},
\text{ for }\|u-v_2\|<\frac{1}{m}.
\]
Taking $u=v_2-h$, we have
\[
\bigg\|\frac{f(x+sv_1+sv_2)-f(x)}{s}-(y+z)\bigg\|<\frac{2}{l}.
\] 
This choice contradicts the choice of $s$.
\end{proof}

Suppose that $X,Y$ are Banach spaces, $f:X\to Y$.
For $x\in X$, $0\neq v\in X$, and $\e>0$ by $O(f,x,v,\e)$ we denote 
the expression 
\[ \sup\bigg\{\bigg\|\frac{f(x+tv)-f(x)}{t}-\frac{f(x+sv)-f(x)}{s}\bigg\|
:0<|t|,|s|<\e\bigg\}.\]
We also define 
\[O(f,x,v):=\lim_{\e\to0+} O(f,x,v,\e).\]
We borrow this definition from~\cite{PZ}.
The following is true in general (in~\cite[Lemma~11]{PZ} it is assumed that $f$ is Lipschitz,
but it is clearly not necessary):
\begin{equation}\label{Ocharder}
f'(x,v)\text{ exists if and only if }O(f,x,v)=0.
\end{equation}

For the rest of this section, $X$ will be a separable Banach space and $Y$
will be a Banach space with RNP. Also, $G\subset X$ will be a closed set
and $f:X\to Y$ a mapping such that there exist $L,\delta>0$ with
\begin{equation}\label{bong3eq}
\|f(y)-f(x)\|\leq L\|y-x\|\quad\text{whenever }y\in G,\ x\in B(y,\delta).
\end{equation}
We also assume that $D$ is a Borel subset of $G$ such that the
distance function $d_G(x)$ is G\^ateaux differentiable at each point $x\in D$.

\begin{lemma}\label{Oborel}
Let $X$ be separable, $0\neq v\in X$,
and we put $g(x):=O(f,x,v)$. Then $g|_{D}$ is Borel measurable.
\end{lemma}

\begin{proof}
Let $w\in D$. Then $h=f|_{B(w,\delta/4)\cap G}$ is $L$-Lipschitz by~\eqref{bong3eq},
and thus $Z=h(B(w,\delta/4)\cap G)$ is separable. Thus, $Z$ can be isometrically embedded
into $\ell_\infty$, and by~\cite[Lemma~1.1(ii)]{BL},
$h$ can be extended to an $L$-Lipschitz mapping $H:X\to\ell_\infty$
(we identify $Z$ with its isometric representation in $\ell_\infty$ for the moment).
By~\cite[Lemma~11(ii)]{PZ}, $G(x):=O(H,x,v)$ is a~Borel measurable function on~$X$.
We will prove that $g(x)=G(x)$ for all $x\in B(w,\delta/4)\cap D$, and conclude 
that $g|_D$ is Borel measurable (by separability of $X$).
\par
Let $x\in B(w,\delta/4)\cap D$. Fix $\gamma>0$ such that $B(x,2\gamma)\subset B(w,\delta/4)$.
Let $\e>0$ and find $0<\tau<\e$ such
that $d_G(x+tv)<\frac{\e}{L}|t|$ and $x+tv\in B(x,\gamma$) whenever $0<|t|<\tau$. Take $\eta:=\frac{1}{2}\min\big(\e,\tau,\frac{L\gamma}{\e}\big)$.
For $0<|s|,|t|<\eta$ find $y,z\in G\cap B(w,\delta/4)$ such that $\|x+tv-y\|<\frac{\e}{L}|t|$ and
$\|x+sv-z\|<\frac{\e}{L}|s|$. Then we have
\[\bigg\|\frac{f(x+tv)-f(y)}{t}\bigg\|\leq\frac{L}{|t|}\|x+tv-y\|\leq\e,\]
and similarly
$\big\|\frac{f(x+sv)-f(z)}{s}\big\|\leq\e$.
Also,
\[\bigg\|\frac{H(y)-H(x+tv)}{t}\bigg\|\leq \frac{L}{|t|}\|x+tv-y\|
\leq \e,\]
and $\big\|\frac{H(x+sv)-H(z)}{s}\big\|\leq\e$. 
Thus using $f(x)=H(x)$, $f(y)=H(y)$, and $f(z)=H(z)$, we obtain
\begin{multline}\label{Hfodh} 
\bigg\|\frac{H(x+tv)-H(x)}{t}-\frac{H(x+sv)-H(x)}{s}\bigg\|\\
\leq \bigg\|\frac{f(x+tv)-f(x)}{t}-\frac{f(x+sv)-f(x)}{s}\bigg\|+\bigg\|\frac{f(x+tv)-f(y)}{t}\bigg\|\\
+\bigg\|\frac{f(x+sv)-f(z)}{s}\bigg\|
+\bigg\|\frac{H(y)-H(x+tv)}{t}\bigg\|+\bigg\|\frac{H(x+sv)-H(z)}{s}\bigg\|\\
\leq O(f,x,v,\e)+4\e.
\end{multline}
By taking a supremum over $0<|s|,|t|<\eta$ in~\eqref{Hfodh}, we
obtain $O(H,x,v,\eta)\leq O(f,x,v,\e)+4\e$. Send $\eta\to0+$ to
get $O(H,x,v)\leq O(f,x,v,\e)+4\e$, and then $\e\to0+$ to see that
$O(H,x,v)\leq O(f,x,v)$.
\par
By~\eqref{bong3eq} and $H$ being $L$-Lipschitz, we can reverse the r\^oles of $f$ and $H$
in the above argument to show that $O(f,x,v)\leq O(H,x,v)$.
\end{proof}

\begin{lemma}\label{likeOPZlem}
If $x\in D$, $0\neq v\in X$, $ O(f,x,v)>0$, $\vp:\R\to X$, $r\in\R$, $\vp(r)=x$,
and the mapping $\psi:t\to\vp(t)-tv$ has Lipschitz constant strictly less
than $O(f,\vp(r),v)/8L$, then the mapping $f\circ\vp$ is not differentiable
at $r$.
\end{lemma}

\begin{proof}
Denote $K:=O(f,x,v)>0$.
To prove the lemma, let $\delta'>0$ be such that $x+tv\in B(x,\delta/2)$ and
$d_G(x+tv)< \frac{K}{16L} |t|$ for each $0<|t|<\delta'$.
Fix $\e>0$ and let $\tau=\min\big(\e,\delta',\frac{16L\delta}{2K}\big)$. By the assumptions on $f$,
let $0<|t|,|s|<\tau$ such that
\[ \bigg\|\frac{f(x+tv)-f(x)}{t}-\frac{f(x+sv)-f(x)}{s}\bigg\|>\,\frac{3}{4}O(f,x,v),\]
and estimate
\begin{equation*}
\begin{split}
D&:= \bigg\|\frac{f\circ\vp(r+t)-f\circ\vp(r)}{t}-\frac{f\circ\vp(r+s)-f\circ\vp(r)}{s}\bigg\|\\
&\geq\bigg\|\frac{f(x+tv)-f(x)}{t}-\frac{f(x+sv)-f(x)}{s}\bigg\|
-
\bigg\|\frac{f(x+tv)-f(\vp(r+t))}{t}\bigg\|\\
&\quad-\bigg\|\frac{f(x+sv)-f(\vp(r+s))}{s}\bigg\|.
\end{split}
\end{equation*}
Find $y,z\in G\cap B(x,\delta)$ such that $\|x+tv-y\|<\frac{K}{16L}|t|$ and $\|x+sv-z\|<\frac{K}{16L}|s|$. 
Then we have
$\big\|\frac{f(x+tv)-f(y)}{t}\big\|\leq
\frac{L}{|t|}\|x+tv-y\|\leq\frac{K}{16}$,
and similarly
\begin{equation*} 
\begin{split}
\bigg\|\frac{f(y)-f(\vp(r+t))}{t}\bigg\|&
\leq \frac{L}{|t|}\|y-\vp(r+t)\|\\
&\leq\frac{L}{|t|}\|y-(x+tv)\|+\frac{L}{|t|}\|\vp(r)+tv-\vp(r+t)\|\\
&\leq\frac{K}{16}+\frac{L}{|t|}\|\psi(r)-\psi(r+t)\|\\
&\leq\frac{K}{16}+L\Lip(\psi)
<\frac{K}{16}+\frac{K}{8}=\frac{3K}{16}.
\end{split}
\end{equation*}
Thus
\begin{equation*}
\begin{split} 
\bigg \|\frac{f(x+tv)-f(\vp(r+t))}{t}\bigg\|&\leq \bigg\|\frac{f(x+tv)-f(y)}{t}\bigg\|+\bigg\|\frac{f(y)-f(\vp(r+t))}{t}\bigg\|\\
&<\frac{K}{16}+\frac{3K}{16}=\frac{K}{4}.
\end{split}
\end{equation*}
Since an analogous estimate holds for $\big\|\frac{f(x+sv)-f(\vp(r+s))}{s}\big\|$,
we obtain
$D>\frac{3}{4}K-2\frac{K}{4}=\frac{O(f,x,v)}{4}$;
so $O(f\circ\vp,r,1)\geq O(f,\vp(r),v)/4$ is strictly positive as required.
\end{proof}

\begin{lemma}\label{bonglem2}
For each $0\neq u\in X$, the set
$\Delta=\{x\in D:f'(x,u)\text{ does not exist}\}$
belongs to $\tilde\mcA(u)$.
\end{lemma}

\begin{proof} 
Since $\Delta=\{x\in D:O(f,x,u)>0\}$ by~\eqref{Ocharder}, and
by Lemma~\ref{Oborel} we have that 
$g(x)=O(f,x,u)$ is Borel on $D$,
we obtain that $\Delta$ is Borel.
By the same reasoning, each $A_k=\big\{x\in\Delta:O(f,x,u)>\frac{1}{k}\big\}$ is Borel for $k\in\N$,
and we have $\Delta=\bigcup_k A_k$. 
To finish the proof of the lemma, it
is enough to show that $A_k\in\tilde\mcA(u,1/16kL)$ for each $k\in\N$.
\par
Let $k\in\N$ be fixed. If $\vp:\R\to X$ is such that the function
$t\to\vp(t)-tu$ has Lipschitz constant at most $1/16kL$, then
Lemma~\ref{likeOPZlem} implies that $f\circ\vp$ is not differentiable
at any $t$ for which $\vp(t)\in A_k$. Hence $B_k:=\{t\in\R:\vp(t)\in A_k\}$
is a subset of the set of points at which $f\circ\vp$ is not differentiable.
Since $f\circ\vp$ is pointwise Lipschitz at all $t$ such that $\vp(t)\in\Delta$,
and since $Y$ has RNP, \cite[Proposition~1]{B} implies that $\lambda(B_k)=0$
as required for showing that $A_k\in\tilde A(u,1/16kL)$.
\end{proof}

\begin{lemma}\label{bong7lem}
Let $X$ be separable.
Then there exists a set $R\in\tilde\mcA$ such that $(N_f\cap D)\setminus R\in\tilde\mcA$,
where $N_f$ is the set of all points $x\in X$ at which $f$ is not G\^ateaux differentiable.
\end{lemma}

\begin{proof} Let $w\in D$, and denote $D_w=D\cap B(w,\delta/4)$.
If $g:=f|_{B(w,\delta/4)\cap G}$, then 
$g$ is $L$-Lipschitz on its domain (by~\eqref{bong3eq}). 
Since $T:=g(B(w,\delta/4)\cap G)$ is separable, we will show that
\[ Z:=\overline{\lin}\{u\in Y:u=f'(x,v)\text{ for some }x\in D_w,v\in X\setminus\{0\}\}\]
is a subset of $W:=\overline{\lin}(T)$ (and thus is separable).
Suppose that $x\in D_w$, $0\neq v\in X$, and $f'(x,v)$ exists.
Fix $\gamma>0$ such that $B(x,2\gamma)\subset B(w,\delta/4)$.
Let $\e>0$ and find $\tau>0$ such that for $0<|t|<\tau$ we have $d_G(x+tv)<\frac{\e}{L}|t|$,
$x+tv\in B(x,\gamma)$,
and $\big\|\frac{f(x+tv)-f(x)}{t}-f'(x,v)\big\|<\e$.
Let $\eta=\min\big(\tau,\frac{L\gamma}{2\e}\big)$ and $0<|t|<\eta$. Find $y\in G\cap B(w,\delta/4)$ with 
$\|x+tv-y\|<\frac{\e}{L}|t|$. Then
\begin{equation*}
\begin{split} \bigg\| f'(x,v)-\frac{f(y)-f(x)}{t}\bigg\|
&\leq \e+\bigg\| \frac{f(x+tv)-f(x)}{t}-\frac{f(y)-f(x)}{t}\bigg\|\\
&\leq \e+\frac{L}{|t|}\|x+tv-y\|\leq 2\e.
\end{split}
\end{equation*}
Since $\frac{f(y)-f(x)}{t}\in W$, send $\e\to0+$ to obtain $d_W(f'(x,v))=0$,
and thus $f'(x,v)\in W$.
\par
Since $X,Z$ are separable, by $R_w$ denote the set obtained as a union
of all $A(k,l,m,y,y')\cap D$ (see Lemma~\ref{likebwlem}) where $k,l,m\in\N$,
$y,y'$ are chosen from a countable dense subset of $Z$ and $v_1,v_2$ are
chosen from a countable dense subset of $X$. By Lemmas~\ref{likebwlem} and~\ref{borellem},
there exists $R_w'\in\tilde\mcA$ 
such that $R_w\subset R_w'$.
We have the following: if $x\in D_w\setminus R'_w$,
then the following implication holds:
\begin{itemize}
\item[$(*)$] If the directional derivative $f'(x,u)$ exists in all 
directions $u$ from a set $U_x\subset X$ whose linear span is dense in $X$,
then $f'(x,v)$ exists for all $v\in\lin_\Q U_x$\footnote{Here, $\lin_\Q V=\{\sum^n_{i=1}q_i v_i:q_i\in\Q,\ v_i\in V,i=1,\dots,n,\ n\in\N\}$.}; 
furthermore, $f'(x,\cdot)$ is
bounded and linear on $\lin_\Q U_x$.
\end{itemize}
The proof of~$(*)$ is similar to the proof of~\cite[Theorem~2]{PZ} and so we omit it.
\par
For the rest of the proof, let $(v_n)_n$ be a complete sequence in~$X$.
Let $\Delta_n=\Delta_n(w)$ be the set $\Delta$ from Lemma~\ref{bonglem2} applied to $v_n$; the lemma
implies that $\Delta_n$ is Borel and $\Delta_n\in\tilde\mcA(v_n)$ for each $n\in\N$.
Denote $F_w=D_w\setminus(\bigcup_n\Delta_n)$.
It follows that $H_w:=F_w\setminus R_w'$ is Borel.
We will show that $f$ is G\^ateaux differentiable at each $x\in H_w$.
\par
Let $x\in H_w$. Fix $\gamma>0$ such that $B(x,2\gamma)\subset B(w,\delta/4)$. 
Let $Q:=\lin_\Q\{v_n:n\in\N\}$.
By~$(*)$ we have a bounded linear mapping $\hat T:Q\to Z$
such that $\hat T(q)=f'(x,q)$ for each $q\in Q$. By the density of $Q$,  
$\hat T$ extends to a~bounded linear mapping
$T:X\to Y$. We have to show that $f'(x,v)=T(v)$ for each $0\neq v\in X$.
Given $0\neq v\in X$ and $\e>0$, by the density of $Q$ and continuity
of $T$ there exists $q\in Q$ such that 
\begin{equation}\label{bong11eq}
\|v-q\|<\frac{\e}{9L}\text{ and }\|T(v-q)\|<\frac{\e}{3}.
\end{equation}
By the existence of $f'(x,q)$ and by the differentiability of the distance function
$d_G(x)$ at the point $x$, there exists $\tau_\e>0$ such that
\begin{equation}\label{bong12eq}
\bigg\|\frac{f(x+tq)-f(x)}{t}-f'(x,q)\bigg\|<\frac{\e}{3},
\end{equation}
$x+tv\in B(x,\gamma)$, and $d_G(x+tv)<\frac{\e}{9L}|t|$ for each $0<|t|<\tau_\e$.
Let $0<|t|<\min(\tau_\e,9\gamma L/2\e)$ and let $y\in G\cap B(w,\delta/4)$ be such that
$\|x+tv-y\|<\frac{\e}{9L}|t|$. Then $\|x+tq-y\|\leq \frac{2\e}{9L}|t|$.
Thus we have
\begin{multline}\label{bong13eq}
\bigg\|\frac{f(x+tv)-f(x+tq)}{t}\bigg\|
\leq \bigg\|\frac{f(x+tv)-f(y)}{t}\bigg\|+ \bigg\|\frac{f(x+tq)-f(y)}{t}\bigg\|
\leq \frac{\e}{3}.
\end{multline}
Now since $f'(x,q)=T(q)$, by~\eqref{bong11eq},~\eqref{bong12eq},
and~\eqref{bong13eq} it follows that
\begin{equation*}
\begin{split}
\bigg\|\frac{f(x+tv)-f(x)}{t}-T(v)\bigg\|&\leq
\bigg\|\frac{f(x+tq)-f(x)}{t}-f'(x,q)\bigg\|\\
&\quad+\bigg\|\frac{f(x+tv)-f(x+tq)}{t}\bigg\|+\|T(v-q)\|\leq\e,
\end{split}
\end{equation*}
for each $0<|t|<\tau_\e$. This proves that $f'(x,v)$ exists and $f'(x,v)=T(v)$.
Thus $f$ is G\^ateaux differentiable at~$x$.
\par
Since there exist $w_k\in D$ such that $D=\bigcup_k (D\cap B(w_k,\delta/4))$,
let $R=\bigcup_k R'_{w_k}$ we have that $R$ is Borel and
since
\begin{equation}
\label{boreq}
(N_f\cap D)\setminus R=\bigg(\bigcup_k((N_f\cap D_{w_k})\setminus R_{w_k}')\bigg)\setminus R
=\bigg(\bigcup_k\big(D_{w_k}\setminus H_{w_k})\bigg)\setminus R,
\end{equation}
we also obtain that $(N_f\cap D)\setminus R$ is Borel (strictly speaking, the
right hand side of~\eqref{boreq} depends on the complete sequence~$(v_n)$,
but the left hand side does not so $(N_f\cap D)\setminus R$ is indeed Borel
since a complete sequence in~$X$ clearly exists by the separability of~$X$).
\par
Since we have the following simple observation: if $A\in\tilde\mcA(v)$
and $B\subset X$ is Borel, then $A\setminus B\in\tilde\mcA(v)$;
we can conclude that $(N_f\cap D)\setminus R$ is indeed in~$\tilde\mcA$.
\end{proof}

\section{Main theorem}\label{mainsec}

\begin{theorem}\label{mainthm}
Let $X$ be a separable Banach space and let $Y$ be a Banach space
with the RNP. Given $f:X\to Y$, let $S(f)$ be the set of all points $x\in X$
at which $f$ is pointwise Lipschitz. Then there exists a set $E\in\tilde\mcA$ such that 
$f$ is G\^ateaux differentiable at every point of $S(f)\setminus E$.
\end{theorem}

\begin{proof}
We follow the proof from~\cite{B}.
For each $n\in\N$ let $G_n$ be the set of all $x\in X$ such that
$\|f(x+h)-f(x)\|\leq n\|h\|$ whenever $\|h\|<\frac{1}{n}$.
Lemma~\ref{bong1lem} implies that each $G_n$ is closed, and $S(f)=\bigcup_n G_n$.
Since the distance function $d_{G_n}(x)$ is Lipschitz on~$X$, by~\cite[Theorem~12]{PZ}
there exists a Borel set $M_n$ such that $X\setminus M_n\in\tilde\mcA$ and $d_{G_n}(x)$
is G\^ateaux differentiable on~$M_n$. Let $D_n:=G_n\cap M_n$.
Thus, in particular, $G_n\setminus D_n\in\tilde\mcA$.
By $\Omega_n$ denote the set of all points $x\in D_n$ at which $f$ is not
G\^ateaux differentiable. By Lemma~\ref{bong7lem} applied
to $D_n$ we obtain $R_n\in\tilde \mcA$ such that $\Omega_n\setminus R_n\in\tilde\mcA$.
\par
Define $E:=\big(\bigcup_n (\Omega_n\setminus R_n)\cup R_n\big)\cup\big(\bigcup_n\big(G_n\setminus D_n\big)\big)$.
Then $E\in\tilde\mcA$ by the previous paragraph. If $x\in S(f)\setminus E$, then
there exists $n\in\N$ such that $x\in G_n\setminus E$. The condition $x\not\in E$
implies that $x\not\in G_n\setminus D_n$ and $x\not\in\Omega_n$.
Therefore $x\in D_n\setminus\Omega_n$, and hence $f$ is G\^ateaux differentiable
at $x$.
\end{proof}

\begin{corollary}\label{maincor}
Let $X$ be a Banach space with $X^*$ separable, $Y$ be a Banach space with RNP,
$f:X\to Y$ be pointwise Lipschitz outside some set $C\in\tilde C$ $($or even some set
$D$ which is $\Gamma$-null$)$, $g:X\to\R$ be continuous convex. Then
there exists a point $x\in X$ such that $f$ is G\^ateaux differentiable
at~$x$ and $g$ is Fr\'echet differentiable at~$x$.
\end{corollary}

\begin{proof}
Assume that $f$ is pointwise Lipschitz outside some $C\in\tilde C$.
By Theorem~\ref{mainthm}, there exists $A\in\tilde\mcA$ such that
$f$ is G\^ateaux differentiable at each $x\in X\setminus (A\cup C)$.
By~\cite[Corollary~3.11]{LP} there exists a $\Gamma$-null $B\subset X$ such that
$g$ is Fr\'echet differentiable at each $x\in X\setminus B$.
Since $A\cup C$ is $\Gamma$-null by~\cite[Theorem~2.4]{Znew}, we have
that $A\cup B\cup C$ is $\Gamma$-null and thus there exists $x\in X\setminus(A\cup B\cup C)$.
\par
If $f$ is pointwise Lipschitz outside a $\Gamma$-null set $D$, then the proof proceeds similarly.
\end{proof}

\section{Cone monotone functions}\label{monotonesection}

\begin{lemma}\label{dominate}
Let $X$ be a Banach space, $K\subset X$ be a closed convex
cone with $0\neq v\in\interior(K)$, and $f:X\to\R$ be $K$-monotone.
If $\limsup_{t\to0} |t|^{-1}|f(x+tv)-f(x)|<\infty$, then $f$ is
pointwise-Lipschitz at $x$.
\end{lemma}

\begin{proof}
Without any loss of generality, we can assume that $v+B(0,1)\subset K$; 
then the proof is identical to the proof of~\cite[Lemma~2.5]{Dcone} (note that
there we assume that $f$ is G\^ateaux differentiable at~$x$,
but, in fact, we are only using that $f$ satisfies $\limsup_{t\to0} |t|^{-1}|f(x+tv)-f(x)|<\infty$).
\end{proof}

Let $(X,\|\cdot\|)$ be a normed linear space. We say that $\|\cdot\|$ is {\em LUR at $x\in S_X$}
provided $x_n\to x$ whenever $\|x_n\|=1$, and $\|x_n+x\|\to 2$.
For more information about rotundity and renormings, see~\cite{DGZ}.

\begin{lemma}\label{renorm}
Let $X$ be a separable Banach space, $K\subset X$ be a closed
convex cone,  $v\in\interior(K)\cap S_X$. 
Then there exists a norm $\|\cdot\|_1$ on $X$ which
is LUR at $v$, $x^*\in(X,\|\cdot\|_1)^*$ with $x^*(v)=\|v\|_1=\|x^*\|=1$,
and $\alpha\in(0,1)$ such that $K_1:=\{x\in X:\|x\|_1\leq \alpha x^*(x)\}$ is
contained in $K$.
\end{lemma}

\begin{proof}
The conclusion follows from~\cite[Lemma~II.8.1]{DGZ} (see e.g.\ the
proof of~\cite[Proposition~15]{Dcone}).
\end{proof}

\begin{lemma}\label{monlem} 
Let $X$ be a Banach space, $v\in S_X$, $x^*\in X^*$ such
that $\|v\|=\|x^*\|=x^*(v)=1$,  $\alpha\in(0,1)$.
Let $K_{\alpha,x^*}=\{x\in X:\alpha\|x\|\leq x^*(x)\}$.
Then there exists $\e=\e(K,v)\in(0,1)$ such that if $\vp:\R\to X$ is
a mapping such that $\psi:t\to\vp(t)-tv$ has Lipschitz constant
less than~$\e$, then $s<t$ implies $\vp(s)\leq_{K_{\alpha,x^*}}\vp(t)$.
\end{lemma}

\begin{proof}
Since $x^*(v)=1$, for each $\alpha<\alpha'<1$ we have $v\in\interior(K_{\alpha',x^*})$.
Fix $\alpha'\in(\alpha,1)$.
Let $\e:=\min\big(1,\frac{(\alpha'-\alpha)}{2\alpha'(1+\alpha)}\big)$.
Take $s<t$, $s,t\in\R$. Then
\begin{equation}\label{coneq}
\begin{split}
\alpha'\|\vp(t)-\vp(s)\|&\leq\alpha'\|\vp(t)-tv-(\vp(s)-sv)\|+\alpha'|t-s|\|v\|\\
&\leq \alpha'\e|t-s|+|t-s|x^*(v)\\
&= \alpha'\e|t-s|+x^*(tv-\vp(t)-(sv-\vp(s))+x^*(\vp(t)-\vp(s))\\
&\leq \alpha'\e|t-s|+\|tv-\vp(t)-(sv-\vp(s))\|+x^*(\vp(t)-\vp(s))\\
&\leq (1+\alpha')\e|t-s|+x^*(\vp(t)-\vp(s)).
\end{split}
\end{equation}
As in~\eqref{coneq}, we show that
$x^*(tv-\vp(t)-(sv-\vp(s)))\leq \e|t-s|$,
and from this we obtain $|t-s|(x^*(v)-\e)\leq x^*(\vp(t)-\vp(s))$.
Then~\eqref{coneq} implies that
\[ \alpha'\|\vp(t)-\vp(s)\|\leq \bigg(1+\frac{(1+\alpha')\e}{1-\e}\bigg)x^*(\vp(t)-\vp(s)).\]
The choice of $\e$ shows that $\alpha\|\vp(t)-\vp(s)\|\leq x^*(\vp(t)-\vp(s))$,
and therefore $\vp(t)\geq_{K_{\alpha,x^*}}\vp(s)$.
\end{proof}

We prove the following theorem, which improves~\cite[Theorem~9]{BW}:

\begin{theorem}\label{conethm}
Let $X$ be a separable Banach space, $K\subset X$ be a closed convex cone
with $\interior(K)\neq\emptyset$. Suppose that $f:X\to\R$ is $K$-monotone.
Then $f$ is G\^ateaux differentiable on~$X$ except for a set belonging to~$\tilde\mcC$.
\end{theorem}

\begin{remark}
It is not known whether $\tilde C\subset\tilde\mcA$ (see~\cite[p.\ 19]{PZ}).
If it is true, then Theorem~\ref{conethm} holds also with $\tilde\mcA$ instead of $\tilde\mcC$.
\end{remark}

\begin{proof}
Without any loss of generality, we can assume that $f$ 
is $K$-increasing and lower semicontinuous (we can work with $\underline{f}$ instead
by~\cite[Proposition~17 and Proposition~16(iii)]{BW}, where $\underline{f}(x)=\sup_{\delta>0}\inf_{z\in B(x,\delta)}f(z)$
is the l.s.c.\ envelope of $f$). 
By Lemma~\ref{renorm},
we can also assume that the norm on $X$ is LUR at $v\in S_X$ and
$K=K_{\alpha,x^*}=\{x\in X:\|x\|\leq\alpha x^*(x)\}$
for some $x^*\in X^*$ and $\alpha\in(0,1)$ with $\|x^*\|=x^*(v)=1$.
\par
Find $\eta>0$ such that $B(v,\eta)\subset\interior(v/2+K_{\alpha,x^*})$
(such an $\eta$ exists since obviously $v\in\interior(v/2+K_{\alpha,x^*})$).
Let $x\in X$ be such that $\|x\|=1$ and $\beta\|x\|\leq x^*(x)$ for
some $0<\beta<1$. Since
\[ 1+\beta=1+\beta\|x\|\leq x^*(v)+x^*(x)\leq\|x+v\|,\]
and the norm on $X$ is LUR at $v$, there exists $\beta'\in(\alpha,1)$
such that $K_{\beta',x^*}\cap S(0,1)\subset B(v,\eta)\subset v/2+K_{\alpha,x^*}$
and thus
\begin{equation}\label{conicaleq} 
K_{\beta',x^*}\cap S(0,t)\subset B(tv,\eta t)\subset tv/2+K_{\alpha,x^*}
\end{equation}
for each $t>0$.
Put $B:=\big\{x\in X: \limsup_{t\to0}\frac{|f(x+tv)-f(x)|}{|t|}=\infty\big\}$.
Then Lemma~\ref{dominate} shows that $S(f)=X\setminus B$, and 
Lemma~\ref{bong1lem} shows that $B$ is Borel.
We will show that $B\in\tilde\mcA(v)$.
Let $\vp:\R\to X$ be a mappings such that $\psi(t)=\vp(t)-tv$ has 
Lipschitz constant strictly less than $\e>0$, where 
$\e$~is given by application of Lemma~\ref{monlem} to $K_{\beta',x^*}$.
Suppose that $r\in\R$ satisfies $\vp(r)=x\in B$. Without
any loss of generality, we can assume that there exist $t_k\to0+$ such that
$\frac{f(x+t_k v/2)-f(x)}{t_k/2}\geq k$ (otherwise work with $-f(-\cdot)$).
For each $k$, find $r_k\in\R$ such that $\vp(r_k)\in (x+K_{\beta',x^*})\cap S(x,t_k)$.
Such $r_k$ exist since $\vp(r)=x$, $\|\vp(s)\|\to\infty$ as $s\to\infty$, and 
$\vp(u)\in (x+K_{\beta',x^*})$ by the choice of $\e$.
Then~\eqref{conicaleq} implies that $\vp(r_k)\geq_{K_{\alpha,x^*}} x+t_k v/2$, and
thus $f(\vp(r_k))\geq f(x+t_k v/2)$. Now, since $\psi$ is $\e$-Lipschitz, we have
$(1-\e)|r-r_k|\leq \|\vp(r_k)-\vp(r)\|=t_k$, and thus
\[ k\leq \frac{f(x+t_k v/2)-f(x)}{t_k/2}\leq \frac{2}{1-\e}\cdot \frac{f(\vp(r_k))-f(\vp(r))}{r-r_k}.\]
It follows that $f\circ\vp$ is not pointwise Lipschitz at $r$.
By the choice of $\e$ and Lemma~\ref{monlem}, we have that $f\circ\vp$ is monotone; 
thus $\lambda(\{r\in\R:\vp(r)\in B\})=0$ (since monotone
functions from $\R$ to $\R$ are known to be a.e.\ differentiable),
and $B\in\tilde\mcA(v,\e/2)$.
\par
We proved that $B\in\tilde\mcA(v)$.
By Lemma~\ref{dominate} we have that $S(f)=X\setminus B$. By Theorem~\ref{mainthm}, there exists
a set $A\in\tilde\mcA$ such that $f$ is G\^ateaux differentiable at all $x\in X\setminus (A\cup B)$.
In~\cite[Theorem~9]{BW} it is proved that the set $N_f$ of points of G\^ateaux non-differentiability
of $f$ is Borel, and thus we obtain that $N_f\in\tilde C$ (since $N_f\subset A\cup B$).
\end{proof}

Theorem~\ref{conethm} and~\cite[Proposition~16(iv)]{BW} show that:

\begin{corollary}\label{Hadcor}
Let $X$ be a separable Banach space, $K\subset X$ be a closed convex
cone with $\interior(K)\neq\emptyset$. Suppose that $f$ is $K$-monotone.
Then $f$ is Hadamard differentiable outside of a set belonging to $\tilde\mcC$.
\end{corollary}

We also have the following corollary.

\begin{corollary}\label{corgamma}
Let $X$ be a Banach space with $X^*$ separable, $K\subset X$ be a closed
convex cone with $\interior(K)\neq\emptyset$, 
$f:X\to\R$ be $K$-monotone, $g:X\to\R$ be continuous convex. Then
there exists a point $x\in X$ such that $f$ is Hadamard differentiable
at~$x$ and $g$ is Fr\'echet differentiable at~$x$.
\end{corollary}

\begin{proof}
By Corollary~\ref{Hadcor}, there exists $A\in\tilde\mcC$ such that
$f$ is Hadamard differentiable at each $x\in X\setminus A$.
By~\cite[Corollary~3.11]{LP} there exists a $\Gamma$-null $B\subset X$ such that
$g$ is Fr\'echet differentiable at each $x\in X\setminus B$.
Since $A$ is $\Gamma$-null by~\cite[Theorem~2.4]{Znew}, we have
that $A\cup B$ is $\Gamma$-null and thus there exists $x\in X\setminus(A\cup B)$.
\end{proof}

\section*{Acknowledgment} The author would like to thank to Prof. Lud\v{e}k Zaj\'\i\v{c}ek
for a useful discussion about $\Gamma$-null sets.

\end{document}